\documentclass[11pt]{article}
\usepackage{amsfonts}\usepackage{amssymb}
\usepackage{psfig}
\usepackage{fullpage}
\usepackage{fleqn}
\usepackage{comment}
\usepackage{here}
\usepackage{euscript}
\usepackage{caption}

\newtheorem{theorem}{Theorem}
\newtheorem{lemma}[theorem]{Lemma}

\newtheorem{proposition}[theorem]{Proposition}
\newcommand{\bpr}{\noindent{\sc Proof: }}
\newcommand{\epr}{~\hfill{\rule{1.5mm}{2.4mm}}}

\begin{document}
\title{Simulation of conditioned diffusions}
 
\author{\Large Bernard Delyon, Ying Hu\\
IRMAR, Universit\'e Rennes 1\\ Campus de Beaulieu, 35042 Rennes Cedex, 
France\\ Emails: bernard.delyon@univ-rennes1.fr, ying.hu@univ-rennes1.fr
}

\maketitle

\centerline{\bf  Abstract}

\bigskip

In this paper, we propose some algorithms for the simulation of the distribution of certain diffusions
conditioned on terminal point. We prove that the conditional distribution is absolutely continuous
with respect to the distribution of another diffusion which is easy for simulation, and
the formula for the density is given explicitly.

\section{Introduction}

The aim of this paper is to propose algorithms for the simulation of the distribution of  a diffusion 
\begin{eqnarray*}
&& {\mbox d}x_t=b(t,x_t){\mbox d}t+\sigma(t,x_t){\mbox d}w_t,~~~x_0=u,~~~0\le t\le T,
\end{eqnarray*}
conditioned  on  $x_T=v$, where $b$ and $\sigma$ are given functions with appropriate dimensions, and
$w$ is a standard Brownian motion.

From the  point of view of application, 
this allows to do posterior sampling when the diffusion is
observed at instants $\{t_1,\cdots, t_n\}\subset[0,T]$.

Let us recall that in the usual conditioning (see, e.g., \cite{LZ}), the distribution of the diffusion $x$ conditioned on $x_T=v$ is the same as that of another diffusion $y$ satisfying
\begin{eqnarray*}
&& {\mbox d}y_t=\tilde{b}(t,y_t){\mbox d}t+\sigma(t,y_t){\mbox d}w_t,~~~y_0=u,~~~0\le t\le T,
\end{eqnarray*}
where
$$\tilde{b}(t,x)=b(t,x)+[\sigma\sigma^*](t,x)\nabla_x(\log p(t,x;T,v)),$$
and $p(s,u;t,z)$ is the density of $x_t^{s,u}$. However, this is not suitable for simulations because in general, one does not know the transition density $p$.

\medskip

We will  prove that, in certain cases,  the conditional
distribution of the diffusion  is absolutely continuous with respect to the distribution of another diffusion
which is easy for simulation, and we give the explicit formula for the density. 
This leads to an efficient simulation algorithm.
\medskip

Two different cases will be considered:

\medskip
\begin{itemize}
\item[1.] The matrix  $\sigma(t,x)$ depends only on $t$, and $b$ has the form
$b(t,x)=b_0(t)+A(t)x+\sigma(t)b_1(t,x)$.
\item[2.] The matrix  $\sigma(t,x)$ is uniformly  invertible.
\end{itemize}

This  simulation algorithm can be applied to a problem of parameter estimation for a discretely observed diffusion as in \cite{roberts}.

This paper is organized as follows: In Sect. 2, we recall a Girsanov theorem for unbounded drift
which is essential for our simulation algorithm.
In Sect. 3, we consider Case 1, and in Sect. 4, we consider Case 2.

\section{A Girsanov theorem for unbounded drifts} 

This section is devoted to give a slightly generalized Girsanov theorem which will be used in the next section.
We call a measurable function $F(t,x)$ from $\mathbb R_+\times\mathbb R^d$
to $\mathbb R^n$ locally Lipchitz with respect to $x$, if for any $R>0$, there exists a constant $C_R>0$, such that,
for any $(t, x,y)\in \mathbb R_+\times  \mathbb R^d \times \mathbb R^d$ with $|x|\le R, |y|\le R$,
$$|F(t,x)-F(t,y)|\le C_R |x-y|.$$
And on the metric space $C([0,T];\mathbb R^m)$, we define the filtration $\{{\cal F}_t\}_t$ to be the natural filtration
of the coordinate process.

\begin{theorem}\label{girsanov}
Let $b(t,x)$, $h(t,x)$, $\sigma(t,x)$ be measurable functions from $\mathbb R_+\times\mathbb R^d$
to $\mathbb R^d$,  $\mathbb R^m$,  and $\mathbb R^{d\times m}$ which are locally Lipschitz with respect to $x$;  
consider the following stochastic differential equations:

\begin{eqnarray}
&& \mathrm{d}x_t=b(t,x_t)\mathrm{d}t+\sigma(t,x_t)\mathrm{d}w_t,\label{eqx}\\
&& \mathrm{d}y_t=(b(t,y_t)+\sigma(t,y_t)h(t,y_t))\mathrm{d}t+\sigma(t,y_t)\mathrm{d}w_t,~~y_0=x_0,\label{eqy}
\end{eqnarray}
on the finite interval $[0,T]$.
We assume the existence of strong solution for each equation.
We assume in addition that $h$ is bounded on compact sets.
Then the Girsanov formula holds: for any non-negative Borel function $f(x,w)$
defined on $C([0,T];\mathbb R^d)\times C([0,T];\mathbb R^m)$, one has
\begin{eqnarray}
&&E_y[f(y,w^h)]=E_x[f(x,w)e^{\int_0^Th^*(t,x_t)\mathrm{d}w_t-\frac{1}{2}\int_0^T|h(t,x_t)|^2\mathrm{d}t}]
,\label{densxy}\\
&&E_x[f(x,w)]=E_y[f(y,w^h)e^{-\int_0^Th^*(t,y_t)\mathrm{d}w_t-\frac{1}{2}\int_0^T|h(t,y_t)|^2\mathrm{d}t}],
\label{densyx}
\end{eqnarray}
where $w_t^h=w_t+\int_0^t h(s,y_s)\mathrm{d}s$, and $h^*$ stands for the transpose of $h$.
\end{theorem}

\bpr We assume first that the positive supermartingale 
\begin{eqnarray*}
&&M_t=\exp\left\{\int_0^th^*(s,x_s){\mbox d}w_s-\frac{1}{2}\int_0^t|h(s,x_s)|^2{\mbox d}s\right\}
\end{eqnarray*}
is a martingale under $P_x$ which will be proved later. 
In this case $\tilde{w}_t=w_t-\int_0^th(s,x_s){\mbox d}s$ is a Brownian motion
under $M_TP_x$, leading to a solution $(x,\tilde{w})$  of (\ref{eqy}):
\begin{eqnarray*}
&& {\mbox d}x_t=b(t,x_t){\mbox d}t+\sigma(t,x_t)h(t,x_t){\mbox d}t
+\sigma(t,x_t){\mbox d}\tilde{w}_t.
\end{eqnarray*}

As $b(t,x),h(t,x),\sigma(t,x)$ are locally Lipschitz with respect to $x$, pathwise uniqueness holds for (\ref{eqx}) and
(\ref{eqy}). The standard Girsanov theorem   implies that (\ref{densxy}) holds.

We prove now that $M_t$ is a martingale.
For any $R>0$, consider the stopping time
\begin{eqnarray*}
&&\tau_R=\inf\{t\ge 0: |x_t|\ge R\}\wedge T.
\end{eqnarray*}
Taking into consideration that $h$ is locally bounded, we have, according to the Girsanov theorem for bounded drift:
\begin{eqnarray*}
&&P_{y|\mathcal F_{\tau_R}}=M_{\tau_R}P_{x|\mathcal F_{\tau_R}}.
\end{eqnarray*}
Hence
\begin{eqnarray*}
&&E_x[M_T]\ge E_x[1_{\tau_R=T}M_T]=E_x[1_{\tau_R=T}M_{\tau_R}]=P_y[\tau_R=T]
\end{eqnarray*}
which converges to 1 as $R\rightarrow\infty$. It implies that $E_x[M_T]=1$, and $M$ is a martingale.

Finally,  (\ref{densyx})  follows in the same way.
\epr

\section{Case when $\sigma$ is independent of $x$}
We assume here that $x_t$ has the specific form
\begin{eqnarray}
\label{prx}
&& \mathrm{d}x_t=(\sigma_t h(t,x_t)+A_t x_t+b_t)\mathrm{d}t+\sigma_t \mathrm{d}w_t,~~~x_0=u,
\end{eqnarray}
where $\sigma_t$ and $A_t$ are time dependent deterministic matrices and $h(t,x)$, 
$b_t$ are vector valued with  appropriate dimension. 

\medskip
For example the 2-dimensional process $(x,y)$ which satisfies the following SDE: 
\begin{eqnarray}
&&\mathrm{d}x_t=y_t\mathrm{d}t \label{pos}\\
&&\mathrm{d}y_t=b(t,x_t,y_t)\mathrm{d}t+\sigma \mathrm{d}w_t \label{vit}
\end{eqnarray}
and which is the noisy version of  $\ddot x_t=b(t,x_t,\dot x_t)$, see, e.g. \cite{arnold}.

We shall prove the following result:
\begin{theorem}
Assume that $A_t, b_t$ and $\sigma_t$ are bounded measurable functions of $t$  with values in  $\mathbb R^{d\times d}$, $\mathbb R^d$
and $\mathbb R^{d\times m}$, respectively.
Assume also that  $h(t,x)$ is locally Lipschitz with respect to $x$ uniformly with respect to $t$ with values in $\mathbb R^m$, and  locally bounded; and the SDE(\ref{prx}) has a strong solution. Moreover, we assume that $\sigma$ admits a measurable left inverse almost everywhere\footnote{This requires essentially that $\sigma^*\sigma$ is almost everywhere $>0$}, denoted by $\sigma^+$; and that $h,A,b$ and $\sigma^+$ are left continuous with respect to $t$. Then,

(i) the covariance matrix $R_{st}$ of the Gaussian process $\xi_t$ corresponding to 
(\ref{prx}) with $h=0$ is given by:

$$R_{st}=P_s\int_0^{\min(s,t)} P_u^{-1} \sigma_u\sigma_u^*P_u^{-*}\mathrm{d}u~P_t^*,$$
where
 $$\frac{\mathrm{d}P_t}{\mathrm{d}t}=A_tP_t,~~~~P_0=Id,$$
and $P_u^{-*}=(P_u^{-1})^*$;

(ii) the distribution of the process
\begin{eqnarray}
&&p_t=\xi_t-R_{tT}R_{TT}^+(\xi_T-v)\label{ptdef}
\end{eqnarray}
is the same as  the distribution of $\xi$ conditioned on $\xi_T=v$
($M^+$ stands for the left pseudo-inverse\footnote{
$M^+=(M^*M)^{-1}M^*$ and the symmetric matrix is inverted by diagonalisation
with $1/0=0$} of $M$).
For any nonnegative measurable function $f$, 
\begin{equation}
E[f(x)|x_0=u,x_T=v]
=C E\left[f(p)e^{\int_0^T  h^*(t,p_t)(\sigma_t^+dp_t-\sigma_t^+(A_t p_t+b_t)\mathrm{d}t)-\frac{1}{2}\int_0^T\|h(t,p_t)\|^2\mathrm{d}t}\right],
\label{simusigma2}\\
\end{equation}
where $C$ is a constant depending on $u,v$ and $T$.
\end{theorem}
\bpr
(i) The formula for $R_{st}$ is classic and comes from
$\xi_t=P_t\int_0^t P_u^{-1} (b_u \mathrm{d}u+\sigma_u \mathrm{d}w_u)+P_t\xi_0$, see e.g. \cite{karatzas}.

\medskip
(ii) Let us first recall that if $(Y,Z)$ is a Gaussian vector, the distribution of $Y$ conditioned on
$Z=z_0$ coincides with the distribution of another Gaussian vector $Y-R_{YZ}R_{ZZ}^+(Z-z_0)$, 
where $R_{ZZ}^+$ is the left pseudo-inverse of $R_{ZZ}$;
its covariance is $R_{YY}-R_{YZ}R_{ZZ}^+R_{ZY}$.
Taking $Y$ as the vector $(\xi_{t_1},\cdots,\xi_{t_k})$, and $Z=\xi_T$,
we observe that, defining the process $p$ by (\ref{ptdef}), $(p_{t_1},\cdots,p_{t_k})$ has the same distribution as that of
$(\xi_{t_1},\cdots,\xi_{t_k})$ conditioned on $\xi_T=v$.
And the covariance of $p_t$ is  $C_{st}=R_{st}-R_{sT}R_{TT}^+R_{Tt}$.

\medskip
Denote by $p^v_t$ the process (\ref{ptdef}); in particular 
for any nonnegative measurable function $\varphi(\cdot)$, 
$E[\varphi(\xi)]=\int E[\varphi(p^v)]\mu_T(\mathrm{d}v)$ where $\mu_T$ is the 
distribution of $\xi_T$.
For any nonnegative measurable functions $f$ and $g$,
\begin{eqnarray}\label{eqrev}
E[f(x)g(x_T)]
&=&E[f(\xi)g(\xi_T)
e^{\int_0^T  h^*(t,\xi_t)\mathrm{d}w_t-\frac{1}{2}\int_0^T\|h(t,\xi_t)\|^2\mathrm{d}t}]\nonumber \\
&=&E[f(\xi)g(\xi_T)
e^{\int_0^T  h^*(t,\xi_t)(\sigma_t^+\mathrm{d}\xi_t-\sigma_t^+(A_t \xi_t+b_t)\mathrm{d}t)-\frac{1}{2}\int_0^T\|h(t,\xi_t)\|^2\mathrm{d}t}].
\end{eqnarray}
Given a sequence of partitions $(\Delta_n)_{n\ge 1}$ of $[0,T]$:
$$\Delta_n=\{t_0^n<t_1^n<\cdots<t_{k_n}^n=T\}$$
with $|\Delta_n|=\max_{0\le i\le k_n-1}(t^n_{i+1}-t_i^n)\rightarrow 0$, and a continuous stochastic process $X$, we define:
$$S_n(X)=\sum_{i=0}^{k_n-1} h^*(t_i^n,X_{t_i^n})\sigma_{t_i^n}^+(X_{t_{i+1}^n}-X_{t_i^n}).$$
Then
$$E[|S_n(\xi)-S_m(\xi)|\wedge 1]=\int_{{\mathbb R}^d} E[|S_n(p^v)-S_m(p^v)|\wedge 1]\mu_T(\mathrm{d}v),$$
which implies that $S_n(p^v)$ converges in probability $P\otimes\mu_T$. Hence, we can define $\int_0^T h^*(t,p_t^v)\sigma_t^+\mathrm{d}p_t^v$
as the limit (in probability $P\otimes\mu_T$) of the sequence $S_n(p^v)$. Obviously, this limit is independent of the sequence of partitions $(\Delta_n)_n$ which satisfies $|\Delta_n|\rightarrow 0$.

Finally, defining the continuous function $\Theta_N(x)=N\wedge x,\ x\ge 0,$ we have
\begin{eqnarray*}
\lefteqn{E[\Theta_N(f(\xi)g(\xi_T)
e^{S_n(\xi)-\int_0^T h^*(t,\xi_t)\sigma_t^+(A_t \xi_t+b_t)\mathrm{d}t-\frac{1}{2}\int_0^T\|h(t,\xi_t)\|^2\mathrm{d}t})]}\\
&=&\int_{\mathbb R^d} E[\Theta_N(f(p^v)e^{S_n(p^v)-\int_0^T h^*(t,p_t^v)\sigma_t^+(A_t p^v_t+b_t)\mathrm{d}t-\frac{1}{2}\int_0^T\|h(t,p^v_t)\|^2
\mathrm{d}t}g(v))]\mu_T(\mathrm{d}v).
\end{eqnarray*}
Taking the limit first in $n$ and then in $N$, and returning to (\ref{eqrev}), we deduce:
\begin{eqnarray*}
 E[f(x)g(x_T)]
=\int_{\mathbb R^d} E[f(p^v)e^{\int_0^T  h^*(t,p^v_t)(\sigma_t^+\mathrm{d}p_t^v-\sigma_t^+(A_t p^v_t+b_t)\mathrm{d}t)-\frac{1}{2}\int_0^T\|h(t,p^v_t)\|^2
\mathrm{d}t}
]g(v)\mu_T(\mathrm{d}v)
\end{eqnarray*}
which implies (\ref{simusigma2}) 
and $C$ is the value of the  density of $\mu_T$ with respect to the distribution of $x_T$ at $v$. 
\epr


As the Brownian bridge, we have:
\begin{proposition}
Let us assume that $M_t=\int_t^TP_u^{-1}\sigma_u\sigma_u^*P_u^{-*}\mathrm{d}u$ is positive definite for any $t\in [0,T)$.
Then the distribution of the process $p$ is the same  as that of $q$ which is the solution to the following linear SDE
\begin{equation}
\label{edsq}
\mathrm{d} q_t=A_tq_t\mathrm{d}t+b_t\mathrm{d}t+\sigma_t\sigma_t^*P_t^{-*}M_t^{-1}
(P_t^{-1}(E[\xi_t]-q_t)-P_T^{-1}(E[\xi_T]-v))\mathrm{d}t+\sigma_t \mathrm{d}w_t,
\end{equation}
with $q_0=u$.
\end{proposition}
\bpr
The matrix  $Q_t=P_tM_t$ is solution to  
 $\dot{Q}_t=(A_t-\sigma_t\sigma_t^*P_t^{-*}M_t^{-1}P_t^{-1})Q_t$,
implying that the covariance of $q_t$ can be rewritten as follows: for $s<t$,
\begin{eqnarray*}
Q_s\int_0^sQ_u^{-1}\sigma_u\sigma_u^*Q_u^{-*}\mathrm{d}u~Q_t^*
&=& Q_s\int_0^s M_u^{-1}P_u^{-1}\sigma_u\sigma_u^*P_u^{-*}M_u^{-1}\mathrm{d}u~Q_t^* \\
&=& Q_s\left (M_s^{-1}-M_0^{-1}\right)Q_t^* \\
&=& P_sM_s(M_s^{-1}-M_0^{-1})M_tP_t^*\\
&=& P_s(M_0-M_s)(Id-M_0^{-1}(M_0-M_t))P_t^*\\
&=& C_{st}.
\end{eqnarray*}
On the other hand, from (\ref{ptdef}), the expectation $\bar p_t$ 
of the process $p_t$ satisfies 
\begin{eqnarray*}
\frac{\mathrm{d}}{\mathrm{d}t}\bar p_t-A_t\bar p_t-b_t
&=&-\sigma_t\sigma_t^*P_t^{-*}P_T^*R_{TT}^{-1}(E[\xi_T]-v).
\end{eqnarray*}
Elementary algebra shows 
$P_T^*R_{TT}^{-1}=-Q_t^{-1}(R_{tT}R_{TT}^{-1}-P_t P_T^{-1})$, hence
\begin{eqnarray*}
\frac{\mathrm{d}}{\mathrm{d}t}\bar p_t-A_t\bar p_t-b_t
&=&\sigma_t\sigma_t^*P_t^{-*}Q_t^{-1}(R_{tT}R_{TT}^{-1}-P_t P_T^{-1})(E[\xi_T]-v)\\
&=&\sigma_t\sigma_t^*P_t^{-*}Q_t^{-1}(E[\xi_t]-\bar p_t-P_t P_T^{-1}(E[\xi_T]-v))
\end{eqnarray*}
which is the equation satisfied by $E[q_t]$. The conclusion follows by noting that 
both $p$ and $q$ are Gaussian processes.
\epr

\paragraph{Remark.} $M_t$ is positive definite for any $t\in [0,T)$ if and only if
the pair of functions $(A,\sigma)$ is controllable on $[t,T]$ for any $t\in [0,T)$.
See, e.g. \cite{karatzas} for some discussions.

\paragraph{Example.} Consider the 2-dimensional stochastic differential equation
defined by (\ref{pos},\ref{vit}), where $\sigma\not=0$. Let us assume that $b$ is locally Lipschitz with respect to $(x,y)$, and this equation admits
a strong solution (the strong solution exists if there exists a Lyapunov function, see, e.g. \cite{arnold}). Then we have:
\begin{eqnarray}
&&E[f(x,y)|(x_0,y_0)=u,(x_T,y_T)=v]
=CE\left[f(p,q)e^{\sigma^{-2}\int_0^T b(t,p_t,q_t)\mathrm{d}q_t
-\frac{1}{2\sigma^2}\int_0^T b(t,p_t,q_t)^2\mathrm{d}t}\right]
\label{simusigma1}
\end{eqnarray}
where $(p,q)$ is the following  bridge starting from $(p_0,q_0)=u$:

\begin{eqnarray}
&&\left(\begin{array}{c}
p_t\\ q_t
\end{array}\right)
=\left(\begin{array}{c}
z_t\\\dot z_t
\end{array}\right)
-\frac{t}{T^3}
\left(\begin{array}{cc}
t(3T-2t) & -tT(T-t)\\6(T-t)&T(3t-2T)
\end{array}\right)
\left(\begin{array}{c}
z_T-v_1\\ \dot z_T-v_2
\end{array}\right),\label{pont2d}\\
\mbox{with} &&z_t= u_1+tu_2+\sigma\int_0^tw_s\mathrm{d}s,~~~\dot z_t=u_2+\sigma w_t;\nonumber
\end{eqnarray}
or $(p,q)$ can be chosen as:

\begin{eqnarray}
&&\mathrm{d}p_t=q_t\mathrm{d}t, \nonumber\\
&&\mathrm{d}q_t=\left(-6\frac{p_t-v_1}{(T-t)^2}-2\frac{2q_t+v_2}{T-t}\right)\mathrm{d}t+\sigma \mathrm{d}w_t.
\label{pont2d0}
\end{eqnarray}

\newpage

\section{$\sigma$ invertible, general $b$}

\subsection{Bounded drift}

Let us consider the following SDEs:
\begin{eqnarray}
&& \mathrm{d}x_t=b(t,x_t)\mathrm{d}t+\sigma(t,x_t) \mathrm{d}w_t,~~~x_0=u, \label{edsx}\\
&& \mathrm{d}y_t=b(t,y_t)\mathrm{d}t- \frac{y_t-v}{T-t}\mathrm{d}t \label{edsy}
+\sigma(t,y_t) \mathrm{d}w_t,~~~y_0=u.
\end{eqnarray}

{\bf Remark.} If $b=0$, and $\sigma=Id$, then $x$ is a Brownian motion. It is well known (see, e.g. \cite{karatzas}) that the law of the Brownian motion $x$
conditioned on $x_T=v$ is the same as that of the Brownian bridge $y$ satisfying the following SDE:
$$\mathrm{d}y_t=- \frac{y_t-v}{T-t}\mathrm{d}t
+ \mathrm{d}w_t,~~~y_0=u.$$
The form of SDE(\ref{edsy}) is inspired by the above SDE in order to fit the simplest case: the Brownian bridge case.

\medskip

The objective of this section is to prove that 
the distribution of $x$ (solution of (\ref{edsx}))  conditioned on $x_T=v$ is absolutely continuous with respect to $y$ (solution of (\ref{edsy})) with an explicit
density. We shall assume some regularity conditions on $b$ and $\sigma$ here.

{\bf Assumption 4.1} The functions $b(t,x)$ and $\sigma(t,x)$ are $C^{1,2}$ with values in 
$\mathbb R^{d}$ and $\mathbb R^{d\times d}$ respectively; and the functions $b,\sigma$, together 
with their derivatives, are bounded. Moreover, $\sigma$ is invertible with a bounded inverse.

Let $x^{s,u}$ be the solution of (\ref{edsx}) starting at $s\in [0,T]$. Under Assumption 4.1, $x$ is a strong Markov process
with positive transition density. For $(s,u)\in [0,T]$, we denote $p(s,u;t,z)$ to be the density of $x_t^{s,u}$. Then there exist
constants $m,\lambda,M,\Lambda>0$, such that the density function $p(s,u;t,z)$ satisfies Aronson's estimation \cite{aronson} : for $t> s$,
$$m(t-s)^{-\frac{d}{2}}e^{-\frac{\lambda|z-u|^2}{t-s}}\le p(s,u;t,z)\le M(t-s)^{-\frac{d}{2}}e^{-\frac{\Lambda|z-u|^2}{t-s}}.$$

We first study SDE (\ref{edsy}).

\begin{lemma}\label{lemmaestimate}
Let Assumption 4.1 hold. Then the SDE(\ref{edsy}) admits a unique solution on $[0,T)$. Moreover,
$\lim_{t\rightarrow T} y_t=v, a.s.$
and
$|y_t-v|^2\le C (T-t)\log\log [(T-t)^{-1}+e],a.s.$,
where $C$ is a positive random variable.
\end{lemma}

\bpr
The fact that the SDE(\ref{edsy}) admits a unique solution on $[0,T)$ is classic. 
Applying It\^o's formula to $\frac{y_t-v}{T-t}$, we deduce easily the following:
$$
\frac{y_t-v}{T-t}=\frac{u-v}{T}+\int_0^t (T-s)^{-1}b(s,y_s)\mathrm{d}s
     +\int_0^t (T-s)^{-1}\sigma(s,y_s)\mathrm{d}w_s.
$$

For each $i$, 
$\{\Big(\int_0^t (T-s)^{-1}\sigma(s,y_s)\mathrm{d}w_s\Big)_i,t\ge 0\}=\{\sum_{j=1}^d \int_0^t (T-s)^{-1}\sigma_{ij}(s,y_s)\mathrm{d}w_s^j,t\ge 0\}$ is a continuous
local martingale, and its quadratic variation process
$\tau_t=\int_0^t \sum_{j=1}^d (T-s)^{-2}\sigma_{ij}^2(s,y_s)\mathrm{d}s$ satisfies $\tau_t\rightarrow \infty$ as $t\rightarrow T$, and
$\tau_t\le \frac{c}{T-t}$ for a constant $c>0$. Applying Dambis-Dubins-Schwarz's theorem, for each $i$, there exists a standard one-dimensional Brownian motion $B^i$, such that $$\Big(\int_0^t (T-s)^{-1}\sigma(s,y_s)\mathrm{d}w_s\Big)_i=B^i(\tau_t),t\ge 0.$$
Taking into  consideration of the law of the iterated logarithm for the Brownian motion $B^i$, the conclusion follows easily.
\epr

Now we can state the main theorem of this section.

\begin{theorem}\label{bounded}
Let Assumption 4.1 hold. Then
\begin{eqnarray}\label{final}
& &E[f(x)|x_T=v]\\
&=& CE\left[f(y)\exp\left\{
-\int_0^T\frac{2\tilde y_t^* A_t(y_t)b_t(y_t) \mathrm{d}t
+ \tilde y_t^*(\mathrm{d}A_t(y_t))\tilde y_t 
+\sum_{ij}\mathrm{d}\langle A^{ij}_t(y_t),\tilde y^i_t\tilde y^j_t\rangle}{2(T-t)}
\right\}\right]\nonumber
\end{eqnarray}
where $A(t,y)=(\sigma(t,y)^*)^{-1}\sigma(t,y)^{-1}$, $\tilde y_t=y_t-v$, and $\langle\cdot,\cdot\rangle$ is the quadratic variation of semimartingales.
\end{theorem}

{\bf Remark.} From Lemma 4, the integral in (\ref{final}) is well defined.

\bigskip

\bpr
Let $f(x)$ be an $\mathcal F_t$-measurable nonnegative  function, $t<T$, then
\begin{eqnarray}\label{eqfy}
E[f(y)]&=&
E\left[f(x)\exp\left\{-\int_0^t\frac{(\sigma^{-1}_s(x_s)(x_s-v))^*}{T-s}\mathrm{d}w_s
-\frac{1}{2}\int_0^t\|\frac{\sigma^{-1}_s(x_s)(x_s-v)}{T-s}\|^2\mathrm{d}s\right\}\right]
\end{eqnarray}
On the other hand, It\^o's formula gives:
\begin{eqnarray*}
\lefteqn{
\mathrm{d}~\frac{\|\sigma^{-1}(t,x_t)(x_t-v)\|^2}{T-t}=
2\frac{(x_t-v)^* A(t,x_t)\mathrm{d}x_t}{T-t}
+\frac{\|\sigma^{-1}(t,x_t)(x_t-v)\|^2}{(T-t)^2}\mathrm{d}t
+\frac{d\cdot \mathrm{d}t}{T-t}}\\
&&+\frac{( x_t-v)^*(\mathrm{d}A(t,x_t))(x_t-v)}{T-t}
+\frac{\sum_{ij}\mathrm{d} \langle A^{ij}(t,x_t),(x^i_t-v^i)(x^j_t-v^j)\rangle}{T-t}
\end{eqnarray*}
Combining the above equation with (\ref{eqfy}), we deduce that,
\begin{eqnarray*}
E[f(y)]&=&C C_t E\left[f(x)\exp\left\{
-\frac{\|\sigma^{-1}_t(x_t)(x_t-v)\|^2}{2(T-t)}
+\int_0^t\frac{( x_s-v)^*A_s(x_s)b_s(x_s) }{T-s}\mathrm{d}s
\right.\right.\\&&\left.\left.
+\frac{1}{2}\int_0^t\frac{( x_s-v)^*(\mathrm{d}A_s(x_s))(x_s-v)}{T-s}
+\frac{\sum_{ij}\mathrm{d}\langle A^{ij}_s(x_s),(x^i_s-v^i)(x^j_s-v^j)\rangle}{T-s}
\right\}\right], 
\end{eqnarray*}
where $C>0$ is a constant, and $C_t=(T-t)^{-\frac{d}{2}}$.

Or equivalently, 
\begin{equation}\label{eq:version1}
E[f(y)\varphi_t]
= C C_t E\left[f(x)\exp\left\{
-\frac{\|\sigma(t,x_t)^{-1}(x_t-v)\|^2}{2(T-t)}\right\}\right],
\end{equation}
where
\begin{equation}\label{eq:version2}
\varphi_{t}=\exp\left\{
-\int_0^t\frac{\tilde y^*_s A_s(y_s)b_s(y_s) }{T-s}\mathrm{d}s
-\frac{1}{2}\int_0^t\frac{ \tilde y_s^*(\mathrm{d}A_s(y_s))\tilde y_s }{T-s}
+\frac{\sum_{ij}\mathrm{d}\langle A^{ij}_s(y_s),\tilde y^i_s\tilde y^j_s\rangle}{T-s}
\right\}.
\end{equation} 
Note that $\{\varphi_t,t\in [0,T]\}$ is a well defined continuous process, thanks to Lemma \ref{lemmaestimate}.

Putting $f=1$ in (\ref{eq:version1}), we deduce then:
\begin{equation}\label{eqf}
\frac{\mathbb E[f(y)\varphi_{t}]}{\mathbb E[\varphi_{t}]}=\frac{E\left[f(x)\exp\left\{
-\frac{\|\sigma(t,x_t)^{-1}(x_t-v)\|^2}{2(T-t)}\right\}\right]}{E\left[\exp\left\{
-\frac{\|\sigma(t,x_t)^{-1}(x_t-v)\|^2}{2(T-t)}\right\}\right]}.
\end{equation}

Assuming that $f(x)$ takes the form $f(x)=g(x_{t_1},\cdots,x_{t_N})$, $0<t_1<t_2<\cdots<t_N<T$, $g\in C_b({\mathbb R}^{Nd})$,
and letting $t\rightarrow T$,
from the Lemmas \ref{lemma:ap1} and \ref{lemma:ap2} in the Appendix, we get:

$$\frac{\mathbb E[f(y)\varphi_{T}]}{\mathbb E[\varphi_{T}]}=\mathbb E[f(x)|x_T=v].$$

This completes the proof of the theorem.
\epr

\medskip
{\bf Remark.} For practical implementation, it is useful to note that the second and third terms of the integral in (\ref{final}) are the limit of $\sum \tilde y_{t_k}^*(A(t_k,y_{t_k})-A(t_{k-1},y_{t_{k-1}}))\tilde y_{t_k} \frac{1}{2(T-t_k)}$.

\subsection{Unbounded drift}

Let us now consider the following SDE:

\begin{equation}
\mathrm{d}x_t=b(t,x_t)\mathrm{d}t+\sigma(t,x_t) \mathrm{d}w_t,~~~x_0=u, \label{edsxunbounded}
\end{equation}
where the drift $b$ can be unbounded. We assume instead

{\bf Assumption 4.2} The function $\sigma(t,x)$ is $C^{1,2}$ with values in
 $\mathbb R^{d\times d}$; the function $\sigma$ together 
with its derivatives are bounded; and  $\sigma$ is invertible with a bounded inverse.
The function $b$ is locally Lipschitz with respect to $x$ and is locally bounded. Moreover, the
SDE (\ref{edsxunbounded}) admits a strong solution.

Combining the Theorems \ref{girsanov} and \ref{bounded}, we are able to prove the following
\begin{theorem}\label{thmunbounded}
Let Assumption 4.2 hold, and $y$ be the solution of
\begin{equation}
\mathrm{d}y_t=- \frac{y_t-v}{T-t}\mathrm{d}t 
+\sigma(t,y_t) \mathrm{d}w_t,~~~y_0=u.\label{edsyunbounded}
\end{equation}

Then,
\begin{eqnarray*}
& &E[f(x)|x_T=v]\\
&=& CE\left[f(y)\exp\left\{
-\int_0^T\frac{\tilde y_t^*(\mathrm{d}A(t,y_t))\tilde y_t 
+\sum_{ij}\mathrm{d}\langle A^{ij}(t,y_t),\tilde y^i_t\tilde y^j_t\rangle }{2(T-t)}\right.\right.\\
& &\left.\left.+\int_0^T (b^*A)(t,y_t)\mathrm{d}y_t-\frac{1}{2}\int_0^T|\sigma^{-1}b|^2(t,y_t)\mathrm{d}t\right\}\right],
\end{eqnarray*}
where $A(t,y)=\sigma(t,y)^{-*}\sigma(t,y)^{-1}$, $\tilde y_t=y_t-v$, and $\langle\cdot,\cdot\rangle$ is the quadratic variation of semimartingales.
\end{theorem}
\bpr
Let $\bar{x}$ be the solution of:
\begin{equation}\label{edsxbar}
\mathrm{d}\bar{x}_t=\sigma(t,\bar{x}_t)\mathrm{d}w_t, ~~~\bar{x}_0=u.
\end{equation}

Then, from Theorem \ref{girsanov}, for nonnegative measurable functions $f$ and $g$,
\begin{eqnarray*}
\mathbb E[f(x)g(x_T)]&=&\mathbb E[f(\bar{x})g(\bar{x}_T)e^{\int_0^T (b^* A)(t,\bar{x}_t)\mathrm{d}\bar{x}_t-\frac{1}{2}\int_0^T|\sigma^{-1}b|^2(t,\bar{x}_t)\mathrm{d}t}]\\
&=&\int_{\mathbb R^d}E[f(\bar{x})e^{\int_0^T (b^* A)(t,\bar{x}_t)\mathrm{d}\bar{x}_t-\frac{1}{2}\int_0^T|\sigma^{-1}b|^2(t,\bar{x}_t)\mathrm{d}t}|\bar{x}_T=v]g(v)\mathrm{d}v.
\end{eqnarray*}

It remains to apply Theorem \ref{bounded}.
\epr

{\bf Remark.} If the drift $b$ is bounded, both formulas in Theorems \ref{bounded} and \ref{thmunbounded} are available. Unfortunately, it is
difficult to compare the efficiency of simulation when applying these two formulas.

\newpage
\section{Appendix}

\begin{lemma}\label{lemma:ap1}
Let $0<t_1<t_2<\cdots<t_N<T$, and $g\in C_b({\mathbb R}^{Nd})$. Then, putting
$$\psi_t=\exp\left\{
-\frac{\|\sigma(t,x_t)^{-1}(x_t-v)\|^2}{2(T-t)}\right\},$$
\begin{equation}\label{eq:ap1}
  \lim_{t\rightarrow T} \frac{\mathbb E[g(x_{t_1},x_{t_2},\cdots,x_{t_N})\psi_{t}]}{\mathbb E[\psi_t]}=\mathbb E[g(x_{t_1},x_{t_2},\cdots,x_{t_N})|x_T=v].
\end{equation}
\end{lemma}

\bpr 
For any $t\in (t_N,T)$,
$$
\frac{\mathbb E[g(x_{t_1},x_{t_2},\cdots,x_{t_N})\psi_{t}]}{\mathbb E[\psi_t]}
=\frac{\int_{\mathbb R^d} \Phi_g(t,z)
\exp\left\{
-\frac{\|\sigma(t,z)^{-1}(z-v)\|^2}{2(T-t)}\right\}\mathrm{d}z}
{\int_{\mathbb R^d}\Phi_1(t,z) \exp\left\{
-\frac{\|\sigma(t,z)^{-1}(z-v)\|^2}{2(T-t)}\right\}\mathrm{d}z},
$$
where
$$\Phi_g(t,z)=\int_{\mathbb R^{Nd}} g(z_1,\cdots,z_N)p(0,u;t_1,z_1)\cdots p(t_N,z_N;t,z)\mathrm{d}z_1\cdots \mathrm{d}z_N,$$
which is continuous thanks to Aronson's estimation. Evidently, $\Phi_1(t,z)=p(0,u;t,z)$.

Moreover, applying a simple change of variable $z=v+(T-t)^{\frac{1}{2}}z'$,
\begin{eqnarray*}
& &(T-t)^{-\frac{d}{2}}\int_{\mathbb R^d} \Phi_g(t,z)
\exp\left\{
-\frac{\|\sigma(t,z)^{-1}(z-v)\|^2}{2(T-t)}\right\}\mathrm{d}z\\
&=&\int_{\mathbb R^d} \Phi_g(t,v+(T-t)^{\frac{1}{2}}z')
\exp\left\{
-\frac{\|\sigma(t,v+(T-t)^{\frac{1}{2}}z')^{-1}z'\|^2}{2}\right\}\mathrm{d}z'\\
&\rightarrow&\Phi_g(T,v)\int_{\mathbb R^d} 
\exp\left\{
-\frac{\|\sigma(T,v)^{-1}z'\|^2}{2}\right\}\mathrm{d}z'.
\end{eqnarray*}

Hence,
$$\lim_{t\rightarrow T} \frac{\mathbb E[g(x_{t_1},x_{t_2},\cdots,x_{t_N})\psi_{t}]}{\mathbb E[\psi_t]}=\frac{\Phi_g(T,v)}{\Phi_1(T,v)},$$
from which we deduce (\ref{eq:ap1}) by the Bayes formula, since
$$\Phi_g(T,v)=\int_{\mathbb R^{Nd}} g(z_1,\cdots,z_N)q(z_1,\cdots,z_N,v)\mathrm{d}z_1\cdots \mathrm{d}z_N,$$
where $q$ is the density of $(x_{t_1},\cdots,x_{t_N},x_T)$.
\epr

\begin{lemma}\label{lemma:ap2}
\begin{eqnarray*}
\lim_{t\rightarrow T}E[|\varphi_t-\varphi_T|]=0.
\end{eqnarray*}
\end{lemma}

We need the following two propositions to prove this lemma.

\begin{proposition}\label{prop:ap1}
(i) There exist two constants $c_1>0,c_2>0$, such that
$$c_1\le C_tE[\psi_t]\le c_2, \forall t\in [0,T),$$
where
$$C_t=(T-t)^{-\frac{d}{2}}.$$

(ii) There exists a constant $c_3>0$, such that
$$E[\varphi_t]\le c_3, \ \forall t\in [0,T).$$
\end{proposition}
\bpr
(i) We  note that
$$C_t\mathbb E[\psi_t]=(T-t)^{-\frac{d}{2}}\int \exp\left\{
-\frac{\|\sigma(t,z)^{-1}(z-v)\|^2}{2(T-t)}\right\}p(0,u;t,z)\mathrm{d}z.$$
We get easily the conclusion taking into consideration of Aronson's estimation after a change of variable $z=v+(T-t)^{\frac{1}{2}}z'$.

(ii) It follows from (\ref{eq:version1}) and (i).
\epr

\begin{proposition}\label{prop:ap2}
For any $\varepsilon>0$, there exists an adapted bounded process $\alpha_t$ such that
\begin{eqnarray*}
\mathrm{d}C_t\psi_t=\mathrm{d}M_t+\alpha_t(C_t\psi_t)^{1-\varepsilon}(T-t)^{-h}\mathrm{d}t,
~~~h=\frac{\varepsilon d+1}{2}
\end{eqnarray*}
where $(M_t)_{0\le t<T}$ is a martingale.
\end{proposition}
\bpr Set $\tilde x_t=x_t-v,~p_t=\|\sigma^{-1}(t,x_t)\tilde x_t\|$,  and 
$A_t=\sigma^{-*}(t,x_t)\sigma^{-1}(t,x_t)$. We have
\begin{eqnarray*}
\mathrm{d}\frac{p_t^2}{T-t}&=&2\frac{\tilde x_t^*A_t\mathrm{d}x_t}{T-t}+
\frac{p_t^2}{(T-t)^2}\mathrm{d}t+\frac{d}{T-t}\mathrm{d}t+\frac{\tilde x_t^*(\mathrm{d}A_t)\tilde x_t}{T-t}
+\frac{1}{T-t}\sum_{i,j}\mathrm{d}\langle A^{ij}_t,\tilde x_t^i\tilde x_t^j\rangle\\
&=&2\frac{\tilde x_t^*\sigma(t,x_t)^{-*}\mathrm{d}w_t}{T-t}+
\frac{p_t^2}{(T-t)^2}\mathrm{d}t+\frac{d}{T-t}\mathrm{d}t+r_t\frac{p_t^2+p_t}{T-t}\mathrm{d}t
+\frac{p_t^2}{T-t}r'_t\mathrm{d}w_t,
\end{eqnarray*}
where $r_t$ and $r_t'$ are two adapted bounded processes.
Hence we get:
\begin{eqnarray*}
\mathrm{d}C_t\psi_t&=&\frac{d}{2(T-t)}C_t\psi_t\mathrm{d}t
-\frac{1}{2}C_t\psi_t\mathrm{d}\left(\frac{p_t^2}{T-t}\right)
+\frac{1}{8}C_t\psi_t\mathrm{d}\left\langle\frac{p_t^2}{T-t}\right\rangle\\
&=& \mathrm{d}M_t+C_t\psi_tr''_t\left(\frac{p_t^2+p_t}{T-t}+\frac{p_t^4+p_t^3}{(T-t)^2}\right)\mathrm{d}t,
\end{eqnarray*}
where $r_t''$ is an adapted bounded process. For
 any $\varepsilon>0$, $e^{-\varepsilon\frac{x^2}{2}}|x|^k, k=1,2,3,4,$ are all bounded functions, then  there exists a constant
 $c_\varepsilon>0$ such that
\begin{eqnarray*}
\psi_t^\varepsilon\left(\frac{p_t^2+p_t}{T-t}+\frac{p_t^4+p_t^3}{(T-t)^2}\right)\le 
\frac{c_\varepsilon}{\sqrt{T-t}}.
\end{eqnarray*}
Hence,
\begin{eqnarray*}
\mathrm{d}C_t\psi_t&=&\mathrm{d}M_t+(C_t\psi_t)^{1-\varepsilon}(T-t)^{-h}r'''_tc_\varepsilon \mathrm{d}t,
\end{eqnarray*} 
where $r_t'''$ is still an adapted bounded process.
\epr

\bigskip

Let us now return to the proof of Lemma \ref{lemma:ap2}.

\bpr
First, from Fatou's lemma and Proposition \ref{prop:ap1},
$$E[\varphi_T]\le \liminf_{t\rightarrow T} E[\varphi_t]\le c_3.$$

We choose $t_0\in (0,T)$ which is close enough  to $T$, and $A$ large enough, and put
\begin{eqnarray*}
&&\sigma=\inf\{t_0<t<T, C_t\psi_t\le\frac{1}{A}\}
=\inf\{t_0<t<T, p_t^2\ge 2(T-t)\log\frac{A}{(T-t)^{\frac{d}{2}}}\}.
\end{eqnarray*}
Under the distribution of $x_.$, $\sigma<T$ a.s.  However under the distribution
of $y_.$, $\lim_{A\rightarrow+\infty}\sigma=T$, a.s., taking into consideration of Lemma \ref{lemmaestimate}.
We have, from (\ref{eqf}), 
$$
\frac{E[\varphi_t1_{\sigma<t}]}{E[\varphi_t]}
=\frac{E[\psi_t1_{\sigma<t}]}{E[\psi_t]}
\le \frac{1}{c_1} E[C_t\psi_t1_{\sigma<t}].
$$

On the other hand, from Proposition \ref{prop:ap2} with a fixed  $\varepsilon\in (0,1/d)$,
$$
\mathrm{d}C_t\psi_t
=\mathrm{d}M_t+\alpha_t(C_t\psi_t)^{1-\varepsilon}(T-t)^{-h}\mathrm{d}t,$$
i.e., 
$$C_t\psi_t
=C_\sigma\psi_\sigma+M_t-M_\sigma
+\int_\sigma^t\alpha_s(C_s\psi_s)^{1-\varepsilon}(T-s)^{-h}\mathrm{d}s.$$

Hence,
$$
E[C_t\psi_t1_{\sigma<t}]
\le A^{-1}+\bar\alpha\int_{t_0}^tE[C_s\psi_s1_{\sigma<s}]^{1-\varepsilon}(T-s)^{-h}\mathrm{d}s,\mbox{ with }
~~~~~~\bar\alpha=\sup_t\|\alpha_t\|_\infty.$$

Therefore,  $E[C_t\psi_t1_{\sigma<t}]$ is bounded by $u_t$ which is  the solution of the following differential equation,
\begin{eqnarray*}
\mathrm{d}u_t&=&\bar\alpha u_t^{1-\varepsilon}(T-t)^{-h}\mathrm{d}t,~~~u_{t_0}=A^{-1};
\end{eqnarray*}
and this equation has an explicit solution:
\begin{eqnarray*}
u_t&=&\left\{\frac{\varepsilon\bar\alpha}{1-h}[(T-t_0)^{1-h}-(T-t)^{1-h}]
+A^{-\varepsilon}\right\}^{1/\varepsilon}\le \left\{c_0(T-t_0)^{1-h}+A^{-\varepsilon}\right\}^{1/\varepsilon},
\end{eqnarray*}
where $c_0>0$ is a constant.
We get finally,
\begin{eqnarray*}
\frac{E[\varphi_t1_{t\le \sigma}]}{E[\varphi_t]}
&=&1-\frac{E[\varphi_t1_{\sigma<t}]}{E[\varphi_t]}\\
&\ge &1-\frac{1}{c_1}(c_0(T-t_0)^{1-h}+A^{-\varepsilon})^{1/\varepsilon}.
\end{eqnarray*}
We note that  $\{\varphi_t1_{t\le\sigma}\}_t$ is a  uniformly integrable family due to Novikov's lemma, since
we have
\begin{eqnarray*}
1_{t\le \sigma}\varphi_t\le C\exp\left\{\int_0^{t\wedge \sigma}\frac{|y_s-v|^2}{T-s}v_s\mathrm{d}w_s-\frac{1}{2}\int_0^{t\wedge \sigma}\frac{|y_s-v|^4}{|T-s|^2}|v_s|^2\mathrm{d}s\right\},
\end{eqnarray*}
where for fixed $A$,  $C$ is a positive constant and $v_t$ is an adapted bounded process.

Taking the $\lim\inf_{t\rightarrow T}$, we get,
\begin{eqnarray*}
\frac{E[\varphi_T1_{\sigma=T}]}{\lim\sup_{t\rightarrow T}E[\varphi_t]}
&\ge &1-\frac{1}{c_1}(c_0(T-t_0)^{1-h}+A^{-\varepsilon})^{1/\varepsilon}.
\end{eqnarray*}
Since  $1_{\sigma=T}$ converges to one a.s. as $A\rightarrow\infty$, we get
\begin{eqnarray*}
\frac{E[\varphi_T]}{\lim\sup_{t\rightarrow T}E[\varphi_t]}
&\ge &1-\frac{1}{c_1}(c_0(T-t_0)^{1-h})^{1/\varepsilon}.
\end{eqnarray*}
It remains to let $t_0\rightarrow T$ to get:
$$\limsup_{t\rightarrow T} E[\varphi_t]\le E[\varphi_T].$$

Hence,
$$\lim_{t\rightarrow T} E[\varphi_t]= E[\varphi_T],$$
and we finish the proof by Scheff\'e's lemma (see, e.g. \cite{DM}).

\epr


\newpage 

\begin{thebibliography}{99}
\bibitem{arnold}
{\sc L.~Arnold},
{Random Dynamical Systems}. Springer, Berlin, 1998.

\bibitem{aronson}
{\sc D.G.~Aronson},
{Bounds for the fundamental solution of a parabolic equation}. {\it Bull. Amer. Math. Soc.} {\bf 73} (1967), 890-896.

\bibitem{CR} {\sc G.~Casella, C.P.~Robert},
{Monte Carlo Statistical Methods}. Springer, New York, 1999.

\bibitem{DM} {\sc C.~Dellacherie, P.A.~Meyer},
{Probabilit\'es et Potentiels.} Chapitre I \`a IV. Hermann, Paris, 1975.


\bibitem{Ha}
  {\sc R.Z.~Has'minskii},
  {Stochastic Stability of Differential Equations}. 
  Sijthoff and Noordhoff, 1980.



\bibitem{karatzas}
  {\sc I.\,Karatzas, S.\,Shreve},
  {Brownian Motion and Stochastic Calculus.} Second edition.
  Springer, New York, 1991.


\bibitem{LZ}
{\sc T.J.~Lyons, W.A.~Zheng}, 
On conditional diffusion processes. {\it Proc. Roy. Soc. Edinburgh Sect. A} {\bf 115} (1990), 243-255.





\bibitem{roberts}
{\sc G.O.\,Roberts, O.\,Stramer}, 
On inference for partially observed nonlinear diffusion models using the Metropolis-Hastings algorithms. {\it Biometrika} {\bf 88} (2001), 603-621.
\end{thebibliography}
\end{document}